\documentclass[11pt]{article}
\usepackage{latexsym,a4wide}
\usepackage{amssymb}
\usepackage{amsmath}
\usepackage{color}
\usepackage{graphicx}
\usepackage{amsthm}
\usepackage{indentfirst}
\allowdisplaybreaks%%%

\newcommand \bel {\begin{equation}\label}
\newcommand \ee {\end{equation}}
\newcommand \be {\begin{equation}}

\newcommand \RR {\mathbb R}

\newcommand \LL {\mathbb L}

\newcommand \del \partial

\newcommand \bei {\begin{itemize}}
\newcommand \eei {\end{itemize}}

\newtheorem{theorem}{\color{black}\indent Theorem}[section]

\newtheorem{remark}{\color{black}\indent Remark}[section]

\begin{document}
\large
\title{Finite time blowup explicit solutions for $3$D incompressible Euler equations }
\author{
{\sc Weiping Yan}\thanks{School of Mathematics, Xiamen University, Xiamen 361000, P.R. China. Email: yanwp@xmu.edu.cn.}
\thanks{Laboratoire Jacques-Louis Lions, Sorbonne Université, 4, Place Jussieu, 75252 Paris, France.}
%\newline  
% {\sl Key Words.}  
%%
%{\sl Mathematics Subject Classification}.  
}
\date{25 June 2018}

\maketitle

\begin{abstract}  
This paper constructs a family of explicit self-similar blowup axisymmetric solutions for the $3$D incompressible Euler equations in $\RR^3$.
Those singular solutions admit infinite energy. 
\end{abstract}

%\tableofcontents

%=========================================================================

\section{Introduction and main results} 
\setcounter{equation}{0}

The $3$D incompressible Euler equation describes the motion of ideal incompressible fluid, which takes the following form
\bel{E1-1}
\aligned
&\textbf{v}_t+\textbf{v}\cdot\nabla \textbf{v}=-\nabla P,\\
&\nabla\cdot\textbf{v}=0,
\endaligned
\ee
where $\textbf{v}(t,x):[0,T^*)\times\RR^3\rightarrow\RR^3$ denotes the $3$D velocity field of the fluid, $P(t,x):[0,T^*)\times\RR^3\rightarrow\RR$ stands for the pressure in the fluid, and $T^*$ is a positive constant.
The divergence free condition in second equation of (\ref{E1-1}) guarantees the incompressibility of the fluid. 

The $3$D incompressible Euler equations develop finite time singularity is regarded as one of the most important open problems in mathematical fluid mechanics \cite{C,F}. 
Toward this direction, Beale-Kato-Majda \cite{BKM} established a celebrated criterion of finite time blowup. 
Constantin-Fefferman-Majda \cite{CFM} gave the geometric aspects the 3D Euler flows and asserts that there can be no blowup if the velocity field is uniformly bounded and the vorticity direction is sufficiently well behaved near the point of maximum vorticity. Lou-Hou \cite{LH} showed a numerical computation to $3$D axisymmetric Euler equations, the solutions were observed to develop singularity on the boundary with local self-similar structure. Recently, Hou-Jin-Liu \cite{HJL} derived a  new family of $3$D models from incompressible axisymmetric Euler and Navier-Stokes equations, and studied its singularity property by numerically. 
The local well-posedness of solutions for this equations are widely studied, one can see \cite{EM,K1,K2,Te} for the strong solutions, and \cite{BT1,BT2,B,Le1,Le2, Sc,Sh,W} for the very weak solutions.
We notice that Constantin \cite{Con} found a class of smooth, mean zero initial data for which the solution of $3$D Euler equations becomes infinite in finite time, meanwhile, he gave an explicit formulas of solutions for the $3$D Euler equations by reducing this equations into a local conservative Riccati system in two-dimensional basic square not in $\RR^3$. 
In this paper, we give a family of explicit self-similar blowup axisymmetric solution for the $3$D incompressible Euler equations (\ref{E1-1}). Here the spatial variable $x\in\RR^3$.

We now derive the $3$D incompressible Euler equations with axisymmetric velocity field in the cylindrical coordinate.
Let $\textbf{e}_r$, $\textbf{e}_{\theta}$ and $\textbf{e}_z$ be the cylindrical coordinate
system,
\bel{E1-0}
\aligned
&\textbf{e}_r=({x_1\over r},{x_2\over r},0)^T,\\
&\textbf{e}_{\theta}=({x_2\over r},-{x_1\over r},0)^T,\\
&\textbf{e}_z=(0,0,1)^T,
\endaligned
\ee
where $r=\sqrt{x_1^2+x_2^2}$ and $z=x_3$. Then the $3$D velocity field $\textbf{v}(t,x)$ is called axisymmetric if it can be written as
$$
\textbf{v}=v^r(t,r,z)\textbf{e}_r+v^{\theta}(t,r,z)\textbf{e}_{\theta}+v^z(t,r,z)\textbf{e}_z,
$$
where $v^r$, $v^{\theta}$ and $v^z$ do not depend on the $\theta$ coordinate.

Thus in the cylindrical coordinates, the $3$D Euler equations (\ref{E1-1}) with axisymmetric velocity field can be reduced into a system as follows
\bel{E1-2}
\aligned
&v_t^{\theta}+v^{r}v_r^{\theta}+v^zv_z^{\theta}=-{v^rv^{\theta}\over r},\\
&\omega_t^{\theta}+v^r\omega_r^{\theta}+v^z\omega_z^{\theta}={2\over r}v^{\theta}v_z^{\theta}+{1\over r}v^r\omega^{\theta},\\
&-\Big(\triangle-{1\over r^2}\Big)\phi^{\theta}=\omega^{\theta},
\endaligned
\ee
where the radial and angular velocity fields $v^r$ and $v^{\theta}$ are recovered from $\phi^{\theta}$ based on the Biot-Savart law
\bel{E1-3}
v^r=-\del_z\phi^{\theta},\quad v^z={1\over r}\del_r(r\phi^{\theta}).
\ee

The incompressibility condition becomes
\bel{E1-4}
\del_r(rv^r)+\del_z(rv^z)=0.
\ee

Obviously, there are singular coefficients ${1\over r}$ in system (\ref{E1-2}). The singular point is $r=0$.
It causes many difficulties to solve it directly. But by some observation, we find some useful structure in system (\ref{E1-2}), then by direct computations, we obtain there exist
a family of explicit self-similar blowup solutions for system (\ref{E1-2}). Furthermore, this solutions also gives a family of explicit
blowup axisymmetric solutions for $3$D incompressible Euler equation (\ref{E1-1}). Here is our main result in this paper.

\begin{theorem}
The 3D incompressible Euler equations (\ref{E1-1}) has a family explicit self-similar blowup axisymmetric solution 
\bel{E1-5}
\textbf{v}(t,x)=v^r(t,r,z)\textbf{e}_r+v^{\theta}(t,r,z)\textbf{e}_{\theta}+v^z(t,r,z)\textbf{e}_z,\quad (t,x)\in[0,T^*)\times\RR^3,
\ee
where $T^*>0$ is a positive constant, $\textbf{e}_r,\textbf{e}_{\theta},\textbf{e}_z$ are defined in (\ref{E1-0}), $r=\sqrt{x_1^2+x_2^2}$ and $z=x_3$, and
$$
\aligned
&v^r(t,r,z)={ar\over T^*-t},\\
&v^{\theta}(t,r,z)={kr^{-(1+{1\over a})}\over T^*-t},\\
&v^{z}(t,r,z)={-2az\over T^*-t},\\
\endaligned
$$
where constants $a,k\in\RR/\{0\}$.
\end{theorem}

We supplement the $3$D incompressible Euler equations (\ref{E1-1}) with the initial data
$$
\textbf{v}(0,x)=\textbf{v}_0(x),\quad x\in\RR^3.
$$
The $3$D incompressible Euler equations satisfies the energy identity
$$
\|\textbf{v}(t,x)\|^2_{\LL^2{(\RR^3})}=\|\textbf{v}(0,x)\|^2_{\LL^2{(\RR^3})}.
$$
By directly computation, it follow from (\ref{E1-5}) that
$$
\textbf{v}(t,x)=\Big(v_1(t,x),v_2(t,x),v_3(t,x)\Big)^T,
$$
where
$$
\aligned
&v_1(t,x)={ax_1\over T^*-t}+{k x_2\over r^{2+{1\over a}}(T^*-t)},\\
&v_2(t,x)={ax_2\over T^*-t}-{kx_1\over r^{2+{1\over a}}(T^*-t)},\\
&v_3(t,x)={-2ax_3\over T^*-t}.
\endaligned
$$
This means that
$$
\aligned
&\nabla v_1(t,x)=\Big({a\over T^*-t}-{k(2+{1\over a})x_1x_2\over 2r^{4+{1\over a}}(T^*-t)},~{k(r^2-(2+{1\over a})x_2^2)\over r^{4+{1\over a}}(T^*-t)},~{kx_2x_3(2+{1\over a})\over r^{4+{1\over a}}(T^*-t) }\Big)^T,\\
&\nabla v_2(t,x)=\Big(-{k(r^2-(2+{1\over a})x_1^2)\over r^{4+{1\over a}}(T^*-t)},~{a\over T^*-t}+{k(2+{1\over a})x_1x_2\over 2r^{4+{1\over a}}(T^*-t)},~-{kx_2x_3(2+{1\over a})\over r^{4+{1\over a}}(T^*-t) }\Big)^T,\\
&\nabla v_3(t,x)=\Big(0,~0,~{-2a\over T^*-t}\Big)^T.
\endaligned
$$
Thus let $ i=1,2,3$, we have
$$
|div(v_i(t,x))|_{x=x_0}|\rightarrow\infty\quad as\quad t\rightarrow (T^*)^-,
$$
where $x_0\in\RR^d$ is a fixed point.

Obviously, there is an initial data
$$
\textbf{v}(0,x)=\Big(v_1(0,x),v_2(0,x),v_3(0,x)\Big)^T,
$$
where
$$
\aligned
&v_1(0,x)={ax_1\over T^*}+{k x_2\over (x_1^2+x_2^2)^{1+{1\over 2a}}(T^*)},\\
&v_2(0,x)={ax_2\over T^*}-{kx_1\over (x_1^2+x_2^2)^{1+{1\over 2a}}(T^*)},\\
&v_3(0,x)={-2ax_3\over T^*}.
\endaligned
$$

Hence the finite time of blowup solution takes place in the gradient of velocity field for the $3$D incomressible Euler equations (\ref{E1-1}).
The singularity is type I. The blowup point is $(t,x)=(T^*,x_0)$. We remark that the stability result on this new blowup axisymmetric solutions (\ref{E1-5})
 for the $3$D incomressible Euler equations (\ref{E1-1}) has been studied in \cite{Yan0}.

\begin{remark}
We remark that the $3$D incompressible Navier-Stokes equations in $\RR^3$
$$
\aligned
&\textbf{v}_t+\textbf{v}\cdot\nabla \textbf{v}=-\nabla P+\nu\triangle \textbf{v},\\
&\nabla\cdot\textbf{v}=0,
\endaligned
$$
have two family of explict blowup axisymmetric solutions as follows
$$
\textbf{v}(t,x)=v^r(t,r,z)\textbf{e}_r+v^{\theta}(t,r,z)\textbf{e}_{\theta}+v^z(t,r,z)\textbf{e}_z,\quad (t,x)\in[0,T^*)\times\RR^3,
$$
where $\textbf{e}_r,\textbf{e}_{\theta},\textbf{e}_z$ are defined in (\ref{E1-0}), $r=\sqrt{x_1^2+x_2^2}$ and $z=x_3$, and
$$
\aligned
&v^r(t,r,z)={ar\over T^*-t},\\
&v^{\theta}(t,r,z)={k\over r},\quad or\quad kr(T^*-t)^{2a},\\
&v^z(t,r,z)=-{2az\over T^*-t},
\endaligned
$$
where constants $a,k\in\RR/\{0\}$. 

When $v^{\theta}(t,r,z)=kr(T^*-t)^{2a}$, the solutions can be rewritten as
$$
\textbf{v}(t,x)=\Big(v_1(t,x),v_2(t,x),v_3(t,x)\Big)^T,
$$
where
$$
\aligned
&v_1(t,x):={ax_1\over T^*-t}+kx_2(T^*-t)^{2a},\\
&v_2(t,x):={ax_2\over T^*-t}-kx_1(T^*-t)^{2a},\\
&v_3(t,x):=-{2ax_3\over T^*-t},
\endaligned
$$
with the smooth initial data
$$
\textbf{v}(0,x)=\Big({ax_1\over T^*}+kx_2(T^*)^{2a},~{ax_2\over T^*}-kx_1(T^*)^{2a},~-{2ax_3\over T^*}\Big)^T.
$$
It is well-known that the equation for pressure $P$ is
$$
-\triangle P=\sum_{i,j=1}^3{\del v_i\over\del x_j}{\del v_j\over\del x_i}.
$$
One can see \cite{Yan00} for more details..
\end{remark}

%--------------------------------------------------------------------------------------------------------------------------

\section{Proof of Theorem 1.1}

In this section, we show how to get the explicit self-similar blowup axisymmetric solutions (\ref{E1-5}) for the $3$D incompressible Euler equations (\ref{E1-1}).
The idea of proof is based on the observation on the structure of system (\ref{E1-2})-(\ref{E1-3}) and incompressibility condition (\ref{E1-4}). It is inspired by the work \cite{Yan1,Yan2}.
Let $T^*$ be a positive constant. We now focus on the axisymmetric system (\ref{E1-2})-(\ref{E1-3}) and incompressibility condition (\ref{E1-4}), that is,
\bel{E2-1}
v_t^{\theta}+v^{r}v_r^{\theta}+v^zv_z^{\theta}=-{v^rv^{\theta}\over r},
\ee
\bel{E2-2}
\omega_t^{\theta}+v^r\omega_r^{\theta}+v^z\omega_z^{\theta}={2\over r}v^{\theta}v_z^{\theta}+{1\over r}v^r\omega^{\theta},
\ee
\bel{E2-3}
-\Big(\triangle-{1\over r^2}\Big)\phi^{\theta}=\omega^{\theta},
\ee
\bel{E2-4}
v^r=-\del_z\phi^{\theta},\quad v^z={1\over r}\del_r(r\phi^{\theta}),
\ee
and the incompressibility condition becomes
\bel{E2-5}
\del_r(rv^r)+\del_z(rv^z)=0.
\ee

We observe the structure of incompressibility condition (\ref{E2-5}), then assume that
\bel{E2-6}
\aligned
&v^r(t,r,z)={a r^p\over T^*-t},\\
&v^z(t,r,z)={b z^q\over T^*-t},
\endaligned
\ee
where $a\neq0,b,p,q$ are unknown constants.

Substituting (\ref{E2-6}) into the incompressibility condition (\ref{E2-5}), we get
$$
\del_r({ar^{p+1}\over T^*-t})+\del_z({bz^qr\over T^*-t})=0,
$$
which gives that
$$
a(p+1)r^p+bqrz^{q-1}=0.
$$
so there must be
$$
p=1,\quad q=1,\quad a=-{b\over 2}.
$$
Thus by (\ref{E2-6}), we get 
\bel{E2-7}
\aligned
&v^r(t,r,z)={a r\over T^*-t},\\
&v^z(t,r,z)=-{2a z\over T^*-t},
\endaligned
\ee
where $a$ is the only unknown constant.

We are now to identify function $\phi^{\theta}(t,r,z)$. 
Substituting (\ref{E2-7}) into (\ref{E2-4}), we have
\bel{E2-9}
\aligned
{a r\over T^*-t}&=-\del_z\phi^{\theta},\\
 -{2a z\over T^*-t}&={1\over r}\del_r(r\phi^{\theta}),
\endaligned
\ee
Assume that
\bel{E2-8}
\phi^{\theta}(t,r,z)={\bar{a}rz\over T^*-t},
\ee
where $\bar{a}$ is an unknown constant.

Substituting (\ref{E2-8}) into (\ref{E2-9}), there is
$$
\aligned
{a r\over T^*-t}&=-{\bar{a}r\over T^*-t},\\
 -{2a z\over T^*-t}&={1\over r}\del_r({\bar{a}r^2z\over T^*-t}),
\endaligned
$$
which means that we should require 
$$
a=-\bar{a},
$$
then above two equations holds.

Thus it follows from (\ref{E2-8}) that
\bel{E2-10}
\phi^{\theta}(t,r,z)=-{arz\over T^*-t}.
\ee
Furthermore, by (\ref{E2-3}), we get
\bel{E2-11}
\omega^{\theta}(t,r,z)=0.
\ee

There is only unknown function $v^{\theta}(t,r,z)$. By observation on the structure of $v^r(t,r,z)$, $v^z(t,r,z)$, $\phi^{\theta}(t,r,z)$ and $\omega^{\theta}(t,r,z)$, we are sure the structure of $v^{\theta}(t,r,z)$ should be
\bel{E2-12}
v^{\theta}(t,r,z)={kz^{p_1}r^{q_1}\over T^*-t},
\ee
where $k,p_1,q_1$ are unknown constants.

Substituting (\ref{E2-7}), (\ref{E2-10})-(\ref{E2-12}) into equations (\ref{E2-1})-(\ref{E2-2}), then we get
\bel{E2-13}
{kz^{p_1}r^{q_1}\over (T^*-t)^2}+{aq_1kr^{q_1}z^{p_1}\over (T^*-t)^2}-{2ap_1kz^{p_1}r^{q_1}\over (T^*-t)^2}=-{akz^{p_1}r^{q_1}\over (T^*-t)^2},
\ee
and
\bel{E2-14}
{2\over r}v^{\theta}v_z^{\theta}=0.
\ee

(\ref{E2-14}) means that $v_z^{\theta}=0$, that is,
$$
p_1=0.
$$
By direct computations, we reduce (\ref{E2-13})-(\ref{E2-14}) into 
\bel{E2-15}
a(1+q_1)+1=0,
\ee
%and%
%\bel{E2-16}
%3a^2-a=2k^2p_1z^{2p_1-2}r^{2q_1}-a^2.
%\ee

If equality (\ref{E2-15}) holds, then there must be
\bel{E2-17}
q_1=-(1+{1\over a}).
\ee

Hence we have obtained all the unknown constants. It follows from (\ref{E2-7}), (\ref{E2-10})-(\ref{E2-12}) that
there are explicit self-similar blowup solutions for system (\ref{E2-1})-(\ref{E2-4}) with condition (\ref{E2-5}) as follows
$$
\aligned
&v^r(t,r,z)={ar\over T^*-t},\\
&v^{\theta}(t,r,z)={kr^{-(1+{1\over a})}\over T^*-t},\\
&v^{z}(t,r,z)={-2az\over T^*-t},\\
&\phi^{\theta}(t,r,z)=-{arz\over T^*-t},\\
&\omega^{\theta}(t,r,z)=0,
\endaligned
$$
where constants $a,k\in\RR/\{0\}$.

Furthermore, note that
$$
\textbf{v}(t,x)=v^r(t,r,z)\textbf{e}_r+v^{\theta}(t,r,z)\textbf{e}_{\theta}+v^z(t,r,z)\textbf{e}_z,\quad (t,x)\in[0,T^*)\times\RR^3,
$$
we obtain the 3D incompressible Euler equations (\ref{E1-1}) admits a family of explicit self-similar blowup axisymmetric solution 
$$
\textbf{v}(t,x)=\Big(v_1(t,x),v_2(t,x),v_3(t,x)\Big)^T,
$$
where
$$
\aligned
&v_1(t,x)={ax_1\over T^*-t}+{k x_2\over r^{2+{1\over a}}(T^*-t)},\\
&v_2(t,x)={ax_2\over T^*-t}-{kx_1\over r^{2+{1\over a}}(T^*-t)},\\
&v_3(t,x)={-2ax_3\over T^*-t}.
\endaligned
$$
This completes the proof.

%================================================================================
%================================================================================
%================================================================================
%================================================================================

\

\

\

\textbf{Acknowledgments.} 

The author expresses his sincerely thanks to the BICMR of Peking University and Profes- sor Gang Tian for constant support and encouragement. The author is supported by NSFC No 11771359.

\end{document}